%% file: frontMatter.tex
\documentclass[a4paper,12pt]{article}
\usepackage{sbpo-template}
\usepackage{amsmath,amssymb,amsfonts,amsthm}
\usepackage{url}
\usepackage[square]{natbib}
\usepackage{indentfirst}
\usepackage{fancyhdr}
\usepackage[utf8]{inputenc}
\usepackage{graphicx}
\usepackage{caption}
\usepackage{subcaption}
\usepackage[version=4]{mhchem}
\usepackage{verbatim}
\usepackage{enumitem}
\usepackage{pgf, tikz}
\usepackage[portuguese]{babel}
\usepackage{natbib}
\usepackage{mathtools}
\usepackage{graphicx}
\usepackage{algorithm}
\usepackage{color}
\usepackage[noend]{algpseudocode}
\usepackage{calrsfs}
\usepackage{glossaries}
\usepackage{rotating, lscape, longtable, tabu, booktabs, siunitx, multirow, multicol}
\usepackage{stackengine}

\usepackage{accents}

\newtheorem{theorem}{Theorem}
\newtheorem{lemma}[theorem]{Lemma}

\newtheorem{definition}{Definition}

\begin{document}

\title{A min-max regret approach for the Steiner Tree Problem with Interval Costs} 

\maketitle
\thispagestyle{fancy}

\author{
\name{Iago A. Carvalho}
\institute{Instituto de Computação, Universidade Estadual de Campinas} 
\iaddress{Av. Albert Einstein 1251, Campinas, SP, 13083-852, Brazil}
\email{iagoac@ic.unicamp.br}
}

\author{ 
\name{Amadeu A. Coco}
\institute{Departamento de Computação, Centro Federal de Educação Tecnológica de Minas Gerais} 
\iaddress{Av. Amazonas 7675, Belo Horizonte, MG, 30510-000, Brazil}
\email {amadeuac@cefetmg.br}
}

\author{ 
\name{Thiago F. Noronha}
\institute{Departamento de Ciências da Computação, Universidade Federal de Minas Gerais}
\iaddress{Av. Antônio Carlos 6627, Belo Horizonte, MG,  31270-901, Brazil}
\email {tfn@dcc.ufmg.br}
}

\author{ 
\name{Christophe Duhamel}
\institute{LITIS, Université Le Havre Normandie} 
\iaddress{25 Rue Philippe Lebon, 76600 Le Havre, France}
\email {christophe.duhamel@univ-lehavre.fr}
}

\begin{abstract}
Let $G=(V,E)$ be a connected graph, where $V$ and $E$ represent, respectively, the node-set and the edge-set. Besides, let $Q \subseteq V$ be a set of terminal nodes, and $r \in Q$ be the root node of the graph. Given a weight $c_{ij} \in \mathbb{N}$ associated to each edge $(i,j) \in E$, the Steiner Tree Problem in graphs (STP) consists in finding a minimum-weight subgraph of $G$ that spans all nodes in $Q$. In this paper, we consider the Min-max Regret Steiner Tree Problem with Interval Costs (MMR-STP), a robust variant of STP. In this variant, the weight of the edges are not known in advance, but are assumed to vary in the interval $[l_{ij}, u_{ij}]$. We develop an ILP formulation, an exact algorithm, and three heuristics for this problem. Computational experiments, performed on generalizations of the classical STP instances, evaluate the efficiency and the limits of the proposed methods.
\end{abstract}

\bigskip

\input{main.tex}

\section*{Aknowledgments}
This study was financed in part by the \emph{Coordenação de Aperfeiçoamento de Pessoal de Nível Superior - Brasil} (CAPES) - Finance Code 001, the \emph{Conselho Nacional de Desenvolvimento Científico e Tecnológico - Brasil} (CNPq), the \emph{Fundação de Amparo à Pesquisa do Estado de Minas Gerais - Brasil} (FAPEMIG), and the \emph{Fundação de Amparo à Pesquisa do Estado de São Paulo - Brasil} (FAPESP).

\bibliographystyle{sbpo}
\bibliography{bibsample}

\end{document}

%% file: main.tex
\section{Introduction} \label{sec:intro}
Let $G = (V,E)$ be a connected graph, where $V$ is the set of nodes and $E$ is the set of edges, where each edge $(i, j) \in E$ is associated with a cost $c_{ij} \in \mathbb{N}_+$. Given a set $Q \subset V$ of terminal nodes, a \textit{Steiner tree} is defined as a tree in $G$ that spans all nodes in $Q$ and may contain additional nodes from $V \setminus Q$. The Steiner Tree Problem in graphs (STP)~\citep{Dreyfus1971} consists in finding a minimum cost Steiner tree of $G$.

STP is a well known NP-Hard problem~\citep{Karp1972}. This problem finds practical applications in areas such as telecommunication networks design, computational biology, VLSI design, among others~\citep{Promel2012}. In most of these applications, the cost associated with each edge is not precisely known. In this paper, we investigate how Robust Optimization (RO)~\citep{Kouvelis1997} can be applied to this context.

RO is an approach to deal with uncertain parameters in decision making, where the data variability is represented by deterministic values~\citep{Kasperski2016,Kouvelis1997}. We focus on RO models where the uncertain data is modelled by an interval of possible values. We refer to the book by \cite{Kouvelis1997} for other robust optimization models. In this approach, any realization of a single value for each parameter is considered as a scenario that can happen. The objective is to find a solution that is efficient for all scenarios, usually referred to as a robust solution. The RO criterion used in this work to classify a solution as robust or not is the \textit{min-max regret}. It was proposed by~\cite{Wald1939} in the context of game theory and was adapted to RO by~\cite{Kouvelis1997}. 

In this paper, we introduce a variant of STP where the value of $c_i$ is uncertain. However, it is assumed that this value is in the range $[l_{ij}, u_{ij}]$. This problem is refereed to as the Min-max Regret Steiner Tree Problem with with Interval Costs (MMR-STP) and is defined as follows.

\begin{definition}
A \emph{scenario} $S$ is an assignment of a single value $c_{ij}^S \in [l_{ij}, u_{ij}]$ for each edge $(i, j) \in E$.
\end{definition}

It is worth noting that there are infinitely many scenarios, as $c_{ij}^S$ can assume any real value in $[l_{ij}, u_{ij}]$. Let $\Gamma$ be the set of all scenarios and $\Phi$ be the set of all Steiner trees in $G$.

\begin{definition}
The \emph{cost} of a solution $x \in \Phi$ in a scenario $S \in \Gamma$ is given by
$$
    F(x,S) = \sum_{(i,j) \in E} c_{ij}^S x_{ij}.
$$
\end{definition}

\begin{definition} 
The cost of the optimal SPT solution $z^S$ in the scenario $S$ is denoted by 
$$
  F(z^S, S) = \min\limits_{z \in \Phi} F(z,S) = \min\limits_{z \in \Phi} \sum_{(i, j) \in E} c_{ij}^S z_{ij}.
$$ 
\end{definition}

\begin{definition}[\cite{Kouvelis1997}]
The \emph{regret} of a solution $x \in \Phi$ in a scenario $S \in \Gamma$ is the difference between the cost of $x$ in the scenario $S$ and the cost of $z^S$ in $S$.
\end{definition}

\begin{definition}
The \emph{worst-case} scenario scenario of $x$, i.e., the one where the regret of $x$ is the maximum is denoted by
$$
  S^x = \arg\max\limits_{S \in \Gamma} \left\{F(x,S) - F(z^S, S)\right\}.
$$
\end{definition}

\begin{lemma}[\cite{Averbakh2001}]
Although $|\Gamma|=\infty$, for any min-max regret optimization problem, $S^x$ is such that $c^{S^{x}}_{ij} = u_{ij}$ if $x_{ij} = 1$, and $c^{S^{x}}_{ij} = l_{ij}$ otherwise. That is, $c_i^{S^x} = l_{ij} + (u_{ij} - l_{ij})x_{ij}$, for all $x \in \Phi$ and $(i,j) \in E$.
\end{lemma}

\begin{definition}
The \textit{robust cost} of a solution $x \in \Phi$ is defined as 
$$
  Z(x) = F(x, S^x) - F(z^{S^x}, S^x),
$$
\end{definition}
\noindent i.e. $Z(x)$ is the maximum regret of $x$. It is worth noting that one has to solve an STP in $S^x$ in order to compute $z^{S^x}$. That is, it is NP-Hard to compute the robust cost of a single solution for MMR-STP. 

\begin{definition}
MMR-STP consists in finding the Steiner tree $x^* \in \Phi$ with the smallest robust cost $Z(x^*)$.
\end{definition}

MMR-STP is clearly NP-Hard, as for $l_{ij} = u_{ij} = c_{ij}$ it reduces to SPT. Therefore, in this paper we propose heuristics and exact algorithms for this problem. We evaluate how good are the solutions provided by the state of the art exact and heuristic algorithms designed for min-max regret optimization problems these problems for the case of MMR-SPT.

The remainder of this work is organized as follows. Related works are presented in Section~\ref{sec:related}. Section~\ref{sec:problem} shows the proposed ILP formulation for MMR-STP. Then, an exact algorithm for MMR-STP are described in Section~\ref{sec:algorithms}, while three heuristics are proposed for this same problem in Section~\ref{sec:heuristics}. Computation experiments, which evaluates the proposed exact and heuristic algorithms, are reported in Section~\ref{sec:experiments}. Finally, concluding remarks are drawn in the last section.


\section{Related work} \label{sec:related}

The Steiner Tree Problem in graphs was proposed in~\cite{Dreyfus1971}, and was proven NP-Hard in~\cite{Karp1972}. Several mathematical formulations~\citep{Chopra1994,Goemans1993,Polzin2001}, as well as exact algorithms~\citep{Lucena1998}, and heuristics~\citep{Duin1994,Duin1999}, were proposed and evaluated for this problem. A comparison among several mathematical formulations for STP was shown in~\cite{Polzin2001}. A compendium of STP formulations can be found in~\cite{Goemans1993}. Furthermore, the state-of-the-art algorithms and other recent advances regarding this problem can be found in~\cite{Du2013,Promel2012}.

Many robust counterparts of classical optimization problems have been studied in the literature, such as the Robust Shortest Path Problem~\citep{Karasan2001,Catanzaro2011,Perez2018} and the Robust Minimum Spanning Tree Problem~\citep{Montemanni2006a,Godinho2019}, and the Robust Shortest Path Tree Problem~\citep{Catanzaro2011,Carvalho2016,Carvalho2016b,Carvalho2018}. These problems are NP-hard~\citep{Aissi2009}, despite the fact that their deterministic counterparts can be solved in polynomial time. RO problems whose deterministic counterparts are already NP-hard have also been studied, such as the Robust Restricted Shortest Path Problem~\citep{Assuncao2017}, the Robust Traveling Salesman Problem~\citep{Montemanni2007}, and the Robust Set Covering Problem~\citep{Pereira2013,Coco2015,Coco2016}, and the Robust Knapsack Problem~\citep{Deineko2010,Furini2015}. As is the case of MMR-STP, these problems are particularly harder to solve than other NP-Hard problems, because the complexity of computing the cost of a single solution is at least that of solving the deterministic counterpart, which is itself NP-Hard. 

Some works in the literature also consider robust variations of STP. A Robust Prize-Collecting Steiner Tree Problem in which both edge weights and node prizes are subject to uncertainty was proposed in~\cite{Alvarez2013}. Moreover, a Two-stage Robust Steiner tree was presented in~\cite{Khandekar2008}. In the initial stage, a small subset of terminal nodes is given. In the second stage, the edge weights are increased by a factor $\lambda$ and several scenarios can occur, each one with a new set of terminal nodes. The objective is to minimize the maximum overall cost over all scenarios. 

Other variations of Steiner problems with data uncertainty were handled by means of stochastic programming. The most studied of these problems is the Two-stage Stochastic Steiner Tree in graphs with recourse~\citep{Bomze2010,Fleischer2006,Gupta2005}. In the initial stage, a known probability distribution $\pi$ is set on subsets of nodes and a cost is assigned to each edge of the graph. In the second stage, a subset of nodes materializes, given their prior known distribution, and the cost of each edge is increased by a factor $\lambda$. Then, an additional set of edges can be bought to build a tree that spans all materialized nodes. The objective is to minimize the expected cost of the two-stage solution. A mathematical formulation and an exact algorithm for this problem were presented in~\cite{Bomze2010}. Approximation algorithms were presented in~\cite{Fleischer2006,Gupta2005}.

To the best our knowledge, the min-max regret Robust Steiner Tree Problem in graphs has not been studied in the literature. Therefore, we propose an Integer Linear Programming (ILP) formulation for MMR-SPT, based on STP's bi-directed multi-commodity flow formulation presented in~\cite{Chopra1994}. As this formulation has an exponentially large number of constraints, we extend the Benders-like Decomposition framework of \cite{Montemanni2005} for MMR-SPT. Furthermore, we propose three heuristics based on the framework of \cite{Kasperski2006}: $(i)$ the Algorithm Mean (AM); $(ii)$ the Algorithm Upper (AU); and $(iii)$ the Algorithm Mean Upper. 

\section{An ILP formulation for MMR-STP} \label{sec:problem}

Our ILP formulation for MMR-STP is based on the multi-commodity flow formulation for STP proposed in~\cite{Chopra1994}. Let $G' = (V, A)$ be a directed graph, obtained by bi-directing the edges in $E$. Furthermore, let $r \in Q$ be an arbitrary terminal node, which is referred to as the \emph{root} node. From this data, STP is formulated by means of binary variables $x \in \{0,1\}^{|A|}$, such that $x_{ij} = 1$ if arc $(i,j) \in A$ belongs to the Steiner tree, and $x_{ij} = 0$ otherwise. Besides, we make use of auxiliary binary variables $y  \in \{0,1\}^{|A \times Q|}$, such that $y_{ij}^k = 1$ if arc $(i,j) \in A$ is used to send an unit of flow from the root  $r$ to the terminal $k \in Q$, and $y_{ij}^k = 0$ otherwise. The resulting formulation consists of the objective function~\eqref{eq:STP_objective} and the constraints in~\eqref{eq:STP_flow}--\eqref{eq:STP_yDomain}. 
\begin{align}
    & \min \sum_{(i,j) \in A} c_{ij} x_{ij} & \label{eq:STP_objective} \\
    & s. t. \nonumber \\
    & \sum_{(j,i) \in A} y_{ji}^k - \sum_{(i,j) \in A} y_{ij}^k =
    \left\{
    \begin{tabular}{l}
        1, if  $j = r$\\
        -1, if  $j = k$\\
        0, otherwise
    \end{tabular}
    \right., 
    \forall j \in N, \ k \in Q \label{eq:STP_flow} \\
    & y_{ij}^k + y_{ij}^k \leq x_{ij}, \quad \forall\,(i,j) \in A, k \in Q \label{eq:STP_projection} \\
    & x_{ij} \in \{0,1\}, \quad \forall\,(i,j) \in A \label{eq:STP_xDomain} \\
    & y_{ij}^k \in \{0,1\}, \quad \forall\,(i,j) \in A, k \in Q \label{eq:STP_yDomain}
\end{align}

The objective function~\eqref{eq:STP_objective} minimizes the cost of the arcs in the Steiner tree. The constraints in~\eqref{eq:STP_flow} are the classic flow conservation constraints that enforce a path from the root $r$ to every other terminal $k \in Q$. The inequalities in \eqref{eq:STP_projection} project the variables $y$ into the variables $x$. Besides,  together with~\eqref{eq:STP_flow}, they enforce that $x$ induce a spanning tree of the terminals in $Q$. The domain of the variables $x$ and $y$ are defined by \eqref{eq:STP_xDomain} and \eqref{eq:STP_yDomain}, respectively.

From the STP formulation described above, we have that the polytope that describes the set $\Phi$ of Steiner trees of $G$ can be formulated by ~\eqref{eq:STP_flow}--\eqref{eq:STP_yDomain}. We have that MMR-STP can be written as $$
\min_{x \in \Phi} Z(x) = \min_{x \in \Phi} F(x, S^{x}) - F(z^{S^x}, S^{x}).
$$
Besides, from definition 1.2, we have that
$$
  F(x, S^{x}) = \sum_{(i,j) \in e} u_{ij} x_{ij}, 
$$
and from definition 1.3 and lemma 1.6, we have that
$$
  F(z^{S^x},S^{x}) 
  = \min_{z \in \Phi} \sum_{(i,j) \in A} c^{S^x}_{ij} z_{ij} 
  = \min_{z \in \Phi} \sum_{(i,j) \in A} \big( l_{ij} + (u_{ij} - l_{ij}) x_{ij} \big) z_{ij}.
$$
Therefore, MMR-SPT can be formulated by the 0-1 Bilevel Integer Linear Program defined by the objective function~\eqref{eq:genericObj} and the constraints~\eqref{eq:STP_flow}--\eqref{eq:STP_yDomain}. We note that~\eqref{eq:genericObj} is indeed linear as $x$ is constant in the inner optimization problem.
\begin{equation} \label{eq:genericObj}
\min_{x \in \Phi} \left\{ \sum_{(i,j) \in A} u_{ij} x_{ij} - \min_{z \in \Phi} \sum_{(i,j) \in A} \big( l_{ij} + (u_{ij} - l_{ij} \big)x_{ij})z_{ij} \right\}
\end{equation}

\sloppy
We can then obtain a MIP formulation by linearizing $F(y^{S^{x}}, S^{x})$, as explained in~\cite{Aissi2009}. Let $y_{ij}^k$ and $x_{ij}$ be the binary variables defined in the STP formulation~\eqref{eq:STP_objective}--\eqref{eq:STP_xDomain}. Besides, let variable $\theta \in \mathbb{R}$ be the cost of the Steiner tree in the worst-case scenario defined by variables $z_{ij}$. The resulting formulation is defined by the objective function~\eqref{eq:MMR-STPObjective} and constraints~\eqref{eq:MMR-STPRobustCost}--\eqref{eq:MMR-STPThetaDomain}. 

\fussy
\begin{align}
    & \min \sum_{(i,j) \in A} u_{ij} x_{ij} - \theta & \label{eq:MMR-STPObjective} \\
    & s. t. \nonumber \\
    & \theta \leqslant \sum_{(i,j) \in A} \big( l_{ij} + (u_{ij} - l_{ij} \big)x_{ij})z_{ij},  \quad \forall\,z \in \Gamma \label{eq:MMR-STPRobustCost} \\
    & \sum_{(j,i) \in A} y_{ji}^k - \sum_{(i,j) \in A} y_{ij}^k =
    \left\{
    \begin{tabular}{l}
        1, if  $j = r$\\
        -1, if  $j = k$\\
        0, otherwise
    \end{tabular}
    \right.,~
    \forall\,j \in N, \ k \in Q  \label{eq:MMR-STPFlow} \\
    & y_{ij}^k + y_{ij}^k \leqslant x_{ij}, \quad \forall\,(i,j) \in A, k \in Q \label{eq:MMR-STPyAndZ} \\
    & x_{ij}   \in \{0,1\}, \quad \forall\,(i,j) \in A \label{eq:MMR-STPxDomain} \\
    & y_{ij}^k \in \{0,1\}, \quad \forall\,(i,j) \in A, k \in Q \label{eq:MMR-STPyDomain} \\
    & \theta   \in \mathbb{R} \label{eq:MMR-STPThetaDomain}
\end{align}

The objective function~\eqref{eq:MMR-STPObjective} aims at minimizing the maximum regret. The constraints in~\eqref{eq:MMR-STPFlow}--\eqref{eq:MMR-STPyDomain} are as previously defined for the STP. The inequalities in~\eqref{eq:MMR-STPRobustCost} computes the cost of each solution $z \in \Gamma$ in the worst-case scenario. One can see that we have an inequality~\eqref{eq:MMR-STPRobustCost} for each possible solution. Therefore, the number of these inequalities is exponential. In order to satisfy these inequalities, the value of $\theta$ should not be greater than the cost of any solution $z \in \Gamma$. 
Finally, the constraint in~\eqref{eq:MMR-STPThetaDomain} defines the domain of variable $\theta$.

\section{A Benders-like Decomposition for MMR-STP} \label{sec:algorithms}

The Benders-like Decomposition (Benders) for MMR-STP is inspired by the Benders Decomposition and based on the approaches used to solve other interval data min-max robust optimization problems, as the Robust Set Covering Problem~\cite{Pereira2013} and the Robust Minimal Spanning Tree Problem~\cite{Montemanni2006a}. Benders is based on formulation \eqref{eq:MMR-STPObjective}--\eqref{eq:MMR-STPThetaDomain}. As the number of constraints \eqref{eq:MMR-STPRobustCost} grows exponentially with the number of nodes, they are relaxed in the master problem. At each iteration, one of these constraints is separated and added to the master problem. Benders stops when the lower bound obtained by solving the master problem is equal to the cost of the best (in this case optimal) solution. 

Let $\Gamma^h \subseteq \Gamma$ be the subset of constraints \eqref{eq:MMR-STPRobustCost} that are known in the master problem at iteration $h$ of Benders, and $X^h$ be the optimal solution of this problem. The value of $\theta$ may not be equal to the cost of the optimal solution $Y^h \in \Gamma$ of scenario $s(X^h)$, as constraints \eqref{eq:MMR-STPRobustCost} are relaxed. 
Therefore, a new constraint \eqref{eq:MMR-STPRobustCost}, generated from $Y^h$, must be added to $\Gamma^{h + 1}$ in order to update the value of $\theta$ for $X^h$. $Y^h$ can be obtained by solving a STP subproblem in $s(X^h)$.

In the first iteration, in order to avoid an unbounded master problem, $\Gamma^1$ is initialized with two solutions obtained by the Algorithm Mean and Algorithm Upper heuristics proposed in Section~\ref{sec:heuristics}, as suggested in \cite{Pereira2013}. At each iteration, the master problem and the corresponding STP subproblem (Equations~\eqref{eq:STP_objective}--\eqref{eq:STP_yDomain}) are solved. Given the lower bound $z^h$ obtained by solving the master problem at iteration $h$, if $z^h < min_{l \in \{1, \dots, h\}} \rho^{s(X^h)}(X^h)$, $\Gamma^{h + 1} = \Gamma^{h} \cup \{Y^h\}$ and a new iteration starts. Otherwise, Benders stops since an optimal solution was found.

\section{Heuristics for MMR-STP} \label{sec:heuristics}

A framework for building heuristics that can be applied to any interval data min-max robust optimization problem was introduced in \cite{Kasperski2006}. The complexity of the algorithms developed through this framework are the same of solving the classical counterpart of the robust optimization problem studied. In this work, we applied this framework to develop three heuristics for MMR-STP.

The first heuristic, called Algorithm Mean (AM), uses a \textit{branch-and-bound} algorithm based on the flow formulation presented in \cite{Polzin2001} to solve a STP at the midpoint scenario $s^=$. In $s^=$, the weight of each edge is set to its mean value, \textit{i.e.} $c_{ij}^{s^=} = (u_{ij} + l_{ij})/2$, for all edges $(i, j) \in E$. Next, the \emph{maximum regret} of the computed solution is evaluated and returned. The cost of the solution obtained through AM is bounded by a factor of $2$ from the optimal solution, as proved in \cite{Kasperski2006}.

The second heuristic, called Algorithm Upper (AU), is similar to AM. However, instead of solving a STP for scenario $s^=$, AU solves the STP for the upper scenario $s^+$, where the weight of each edge is set to its upper value, \textit{i.e.} $c_{ij}^{s^+} = u_{ij}$, for all edges $(i, j) \in E$. Unlike AM, the cost of the solution obtained by AU is not bounded.

The last heuristic, called Algorithm Mean Upper (AMU), combines AM and AU. Next, it returns the smallest computed \emph{maximum regret}. As AMU runs AM, it is also a $2$-approximation algorithm for any interval data min-max robust optimization problem.

\section{Computational experiments} \label{sec:experiments}

Computational experiments have been performed on an Intel Xeon CPU E5645 with $2.4$ GHz clock and $32$ GB of RAM memory, running under Linux operating system. The \emph{branch-and-bound} implementation of the ILOG CPLEX version $12.6$ with default parameter settings was used to solve the mixed integer linear programs. The algorithms were implemented in C++ using the ILOG Concert Technology and compiled with GNU g++ 5.4.0. The running time of all algorithms has been limited to 10800 seconds (3 hours). 

The instances used in the experiments are generalizations of classical STP instances. The 5 first instances (WRP3-11 to WRP3-15) from the SteinLib \footnote{http://elib.zib.de/steinlib} WRP3 set were used. Their sizes range from 128 nodes, 227 edges, and 11 terminal nodes (WRP3-11) to 138 nodes, 257 edges, and 15 terminal nodes (WRP3-15). Next, three different methods, namely Beasley (BE), Montemanni (MO) and Kasperski-Zielinski (KZ), are used to generate the edge weights interval as in \cite{Pereira2013}. A parameter $\beta = \{0.1, 0.3, 0.5\}$ is used in BE, while $M = \{750, 1000, 1250\}$ is used as parameter in MO and KZ. For each method, the higher the parameter value, the larger the edge weight interval. These methods are applied to the selected WRP3 instances. Therefore, nine sets of 5 instances were generated by using different interval sizes. They are used in the experiment described below.

The performance of Benders, AM, AU, and AMU for these sets is displayed in Table \ref{tab:results}. The name of each instance set is shown in Column 1. The average relative optimality gap of Benders, as well as the average computation time of these runs are reported in columns 2 and 3, respectively. Then, columns 4 and 5 present respectively $(i)$ the average percent relative deviation to the Bender's upper bound and $(ii)$ the average computational time of AM for each set. The same information is given for AU and for AMU.

\begin{table}[h]
	\centering
	\setlength{\tabcolsep}{.50em}
	\begin{tabular}{lrrrrrrrr} 
	    \toprule
		&       \multicolumn{2}{c}{Benders} & \multicolumn{2}{c}{AM} & \multicolumn{2}{c}{AU} & \multicolumn{2}{c}{AMU} \\ \cmidrule(lr){2-3} \cmidrule(lr){4-5} \cmidrule(lr){6-7} \cmidrule(lr){8-9}
		Instance set & gap\% &    t(s) & \%dev&    t(s) & \%dev&    t(s) & \%dev& t(s)  \\
		\midrule
		WRP3-BE-0.1  & 15.48 & 5948.72 & 3.70 &    1.91 & 1.83 &    1.94 & 1.46 &    3.85 \\
		WRP3-BE-0.3  & 14.80 & 4763.78 & 3.12 &    2.03 & 1.03 &    2.74 & 0.91 &    3.77 \\
		WRP3-BE-0.5  &  6.25 & 4825.14 & 0.85 &    1.76 & 1.65 &    0.97 & 0.14 &    3.26 \\ \midrule
		WRP3-MO-750  &  6.32 & 6059.24 & 4.33 &  238.80 & 2.29 &  240.00 & 1.21 &  478.80 \\
		WRP3-MO-1000 & 12.60 & 7434.67 & 4.22 &  467.28 & 2.80 &  637.26 & 1.96 & 1104.54 \\
		WRP3-MO-1250 & 23.64 & 8583.08 & 5.93 & 1231.83 & 2.79 & 1658.32 & 2.79 & 2890.15 \\ \midrule
		WRP3-KZ-750  & 28.37 &10800.00 & 3.94 & 2916.45 & 3.78 & 2508.25 & 2.02 & 5424.70 \\
		WRP3-KZ-1000 & 27.88 &10800.00 & 1.25 & 2247.66 & 4.10 & 2894.98 & 0.91 & 5142.64 \\
		WRP3-KZ-1250 & 27.92 & 8926.96 & 1.83 & 1510.39 & 3.12 & 1486.42 & 0.93 & 2996.81 \\
		\bottomrule
	\end{tabular}
	\caption{Evaluation of Benders, AM, AU, and AMU}
	\label{tab:results}
\end{table}

One can see from Table~\ref{tab:results} that Benders achieves a smaller relative optimality gap and running times in BE and MO instances than in KZ. It indicates that the KZ instances are the most difficult ones among the three proposed sets. Regarding the interval sizes, we obtained different results for each instance set. For the BE instances, the smaller is the interval size, the greater is the average optimality gap. On the other hand, for the MO instances, the greater is the interval size, the smaller is the average optimality gap. However, the Benders's algorithm average optimality gap was almost the same for all of the KZ instances.

Regarding the heuristics, one can see from this same table that, for BE instances, AM, AU and AMU average running times never exceeds 4 seconds, but grow quickly for MO and KZ instances. The maximum average relative deviations for AM, AU, and AMU are respectively 5.93\%, 4.10\%, and 2.79\%. These results indicate that the Bender's algorithm did not greatly improved its initial solution (which is given by AMU, as explained in Section~\ref{sec:algorithms}).

\section{Conclusions}

This paper considers a new combinatorial optimization problem that arises from the uncertain nature of Steiner tree problem applications. We propose a mathematical formulation based on the robust optimization framework presented in \cite{Kouvelis1997} and an exact and three heuristic algorithms to solve it. Computational experiments show that Benders-like decomposition did not solve all proposed instances to optimality. However, the heuristics achieve good results in a small running time. Future works should focus on the development of new exact and heuristic methods for the studied problem. Moreover, other mathematical formulations for the Steiner tree problem in graphs presented in \cite{Polzin2001} can be extended to handle the uncertain data for this problem.